\documentclass{icm2010}

\title[Multiple orthogonal polynomials]{Multiple orthogonal polynomials
in random matrix theory}
\author[Kuijlaars]{Arno B.J. Kuijlaars
\thanks{The author was supported in part by FWO-Flanders project G.0427.09,
by K.U. Leuven research grant OT/08/33, by the Belgian
Interuniversity Attraction Pole P06/02, and by grant MTM2008-06689-C02-01 of
the Spanish Ministry of Science and Innovation.} }
\contact[arno.kuijlaars@wis.kuleuven.be]{Department of
Mathematics, Katholieke Universiteit Leuven, Celestijnenlaan 200
B, 3001 Leuven, Belgium}

\theoremstyle{definition}

    \DeclareMathOperator*{\diag}{diag}

    \DeclareMathOperator*{\Prob}{Prob}
    \DeclareMathOperator*{\Tr}{Tr}
    \DeclareMathOperator*{\Ai}{Ai}
    \DeclareMathOperator*{\im}{Im}
    \renewcommand{\Im}{\im}

\usepackage{overpic}

\begin{document}

\begin{abstract}
Multiple orthogonal polynomials are a generalization of
orthogonal polynomials in which the orthogonality is distributed
among a number of orthogonality weights. They appear
in random matrix theory in the form of special determinantal point processes
that are called multiple orthogonal polynomial (MOP) ensembles.
The correlation kernel in such an ensemble is expressed
in terms of the solution of a Riemann-Hilbert problem,
that is of size $(r+1) \times (r+1)$ in the case of $r$ weights.

A number of models give rise to a MOP ensemble, and we discuss recent
results on models of non-intersecting Brownian motions, Hermitian
random matrices with external source, and the two matrix model.
A novel feature in the asymptotic analysis of the
latter two models is a vector equilibrium
problem for two or more measures, that describes the limiting
mean eigenvalue density. The vector equilibrium problems involve
both an external field and an upper constraint.
\end{abstract}



\maketitle


\section{Introduction}

\subsection{Random matrix theory}
The Gaussian Unitary Ensemble (GUE) is the most prominent and most studied
ensemble in random matrix theory. It is a probability measure on $n \times n$
Hermitian matrices for which the joint eigenvalue probability density function (p.d.f.)
has the explicit form
\begin{align} \label{GUEdensity}
    \frac{1}{Z_n} \prod_{1\leq j < k \leq n} (x_k-x_j)^2
    \prod_{j=1}^n e^{- \frac{n}{2}  x_j^2}
\end{align}
where $Z_n$ is an explicitly known constant.
The density \eqref{GUEdensity} can be analyzed with the help of
Hermite polynomials. Due to this connection with classical orthogonal
polynomials many explicit calculations can be done,
both for finite $n$ and in the limit $n \to \infty$, see \cite{Meh}.
In particular it leads to a description of the limiting behavior of
eigenvalues on the global (macroscopic) scale
as well as on the local (microscopic) scale. The global scale is given
by the well-known Wigner semi-circle law
\begin{align} \label{Wignersemicircle}
    \rho(x) = \frac{1}{2\pi} \sqrt{4-x^2} \qquad -2 \leq x \leq 2,
\end{align}
in the sense that for eigenvalues $x_1, \ldots, x_n$ taken from \eqref{GUEdensity},
the empirical eigenvalue distribution  $\frac{1}{n} \sum_{j=1}^n \delta(x_j)$
converges weakly to $\rho(x)$
almost surely as $n \to \infty$.

The local scale is characterized by the sine kernel
\begin{equation} \label{sinekernel}
    \mathcal S(x,y) = \frac{\sin \pi(x-y)}{\pi(x-y)}
\end{equation}
in the bulk. This means that for any given $x^* \in (-2,2)$
and any fixed $m \in \mathbb N$, the $m$-point correlation function
(i.e., the marginal distribution)
\begin{multline} \label{correlationfunctions}
    R_{m,n}(x_1, \ldots, x_m) \\
    = \frac{n!}{(n-m)!} \int_{\mathbb R^{n-m}}
     \left[\frac{1}{Z_n}   \prod_{1 \leq j < k \leq n} (x_k-x_j)^2 \prod_{j=1}^n e^{- \frac{1}{2} n x_j^2} \right]
        dx_{m+1} \cdots dx_n \end{multline}
has the scaling limit
\begin{multline} \label{scalinglimit}
    \lim_{n \to \infty}
    \frac{1}{[\rho(x^*) n]^m} R_{m,n} \left(x^* + \frac{x_1}{\rho(x^*)n},
        \ldots, x^* + \frac{x_m}{\rho(x^*)n} \right) \\
        =   \det \left[ \mathcal S(x_i,x_j) \right]_{1 \leq i, j \leq m}.
    \end{multline}

At the edge points $\pm 2$ the sine kernel \eqref{sinekernel}
is replaced by the Airy kernel
\begin{equation} \label{Airykernel}
    \mathcal A(x,y) = \frac{\Ai(x) \Ai'(y) - \Ai'(x) \Ai(y)}{x-y}
\end{equation}
and a scaling limit as in \eqref{scalinglimit} (with scaling factor $c n^{2/3}$
instead of $\rho(x^*) n$) holds for $x^* = \pm 2$.
This result leads in particular to the statement about the largest eigenvalue
\begin{equation} \label{TracyWidom}
    \lim_{n \to \infty} \Prob \left( \max_{1 \leq k \leq n} x_k < 2 +  \frac{t}{cn^{2/3}} \right)
     = \det \left[I - \mathcal A_{(t, \infty)} \right]
    \end{equation}
where $\mathcal A$ is the Airy kernel \eqref{Airykernel}
and the determinant is the Fredholm determinant of the integral
operator with Airy kernel acting on $L^2(t, \infty)$.
The limiting distribution \eqref{TracyWidom} is the famous Tracy-Widom distribution
named after the authors of the seminal work \cite{TW1} in which the right-hand side of \eqref{TracyWidom}
is expressed in terms of the Hastings-McLeod solution of the Painlev\'e II equation.

These basic results of random matrix theory  have
been extended and generalized in numerous directions. Within the theory of
random matrices, they have been
generalized to ensembles with unitary, orthogonal and symplectic symmetry and
to non-invariant ensembles (Wigner ensembles). The distribution functions of
random matrix theory also appear in many other
probabilistic models that have no apparent connection with random matrices
(models of non-intersecting paths, tiling models, and stochastic growth models),
see e.g. \cite{BDJ}, \cite{Joh}.

Mehta's book \cite{Meh} is the standard reference on random matrix theory.
The book of Deift \cite{Dei} has been very influential in introducing
Riemann-Hilbert techniques into the study of random matrices.
In recent years, a number of new monographs appeared \cite{AGZ}, \cite{BS}, \cite{Blo},
\cite{DG}, \cite{For} that cover the various aspects of the theory of
random matrices.

\subsection{Unitary ensembles and orthogonal polynomials}
One direction within random matrix theory is the study
of ensembles of the form
\begin{equation} \label{unitaryensemble}
    \frac{1}{Z_n} e^{-n \Tr V(M)} \, dM
\end{equation}
defined on $n \times n$ Hermitian matrices $M$, which reduces
to the GUE in case $V(x) = \frac{1}{2}x^2$. The ensembles \eqref{unitaryensemble} have the
property of unitary invariance and are called unitary ensembles.
The eigenvalues have the  p.d.f.\
\begin{equation} \label{unitarypdf}
    \frac{1}{Z_n} \prod_{1 \leq j < k \leq n} (x_k-x_j)^2 \prod_{j=1}^n e^{-n V(x_j)}
\end{equation}
with a different normalizing constant $Z_n$. [Throughout, we use $Z_n$ to denote
a normalizing constant, which may be different from one formula to the next.]

Again explicit calculations can be done due to the connection with
orthogonal polynomials \cite{DG}, \cite{Meh}.
For a given $n$, we consider the monic polynomial $P_{k,n}$ of degree $k$
that satisfies
\[ \int_{-\infty}^{\infty} P_{k,n}(x) x^j e^{-n V(x)} dx = h_{k,n} \delta_{j,k}, \qquad
    j=0, \ldots, k. \]
Then  \eqref{unitarypdf} is a determinantal point process \cite{AGZ}, \cite{Sos}, with kernel
\begin{equation} \label{OPkernel}
    K_n(x,y) = \sqrt{e^{-n V(x)}} \sqrt{e^{-n V(y)}}
    \sum_{k=0}^{n-1}  \frac{P_{k,n}(x) P_{k,n}(y)}{h_{k,n}}
    \end{equation}
which means that for every $m \in \mathbb N$ the $m$-point correlation functions,
defined as in \eqref{correlationfunctions},
have the determinantal form
\[ \det \left[ K_n(x_i,x_j) \right]_{i,j=1, \ldots, m}. \]

As $n \to \infty$, the limiting mean eigenvalue density
\[ \rho(x) = \lim_{n \to \infty} \frac{1}{n} K_n(x,x) \]
is no longer Wigner's semi-circle law \eqref{Wignersemicircle},
but instead it is the density $\rho$ of the probability measure $\mu$ that
minimizes the weighted logarithmic energy
\begin{equation} \label{weightedenergy}
    \iint \log \frac{1}{|x-y|} d\mu(x) d\mu(y) + \int V(x) d\mu(x)
\end{equation}
among all probability measures on $\mathbb R$.

Local eigenvalue statistics, however, have a universal behavior as $n \to \infty$,
that is described by the sine kernel \eqref{sinekernel} in the bulk.
Thus for points $x^*$ with $\rho(x^*) > 0$ the limit \eqref{scalinglimit}
holds true.  At edge points of the limiting spectrum the density $\rho$
typically vanishes as a square root and then the universal Airy kernel \eqref{Airykernel}
appears. For real analytic potentials $V$ this was proved in \cite{BI1}, \cite{DKMVZ}
using Riemann-Hilbert methods. This was vastly extended
to non-analytic potentials in recent works of Lubinsky \cite{Lub} and Levin and Lubinsky
\cite{LL}, among many others.

\subsection{This paper}

In this paper we present an overview of the work (mainly of the
author and co-workers) on multiple orthogonal polynomials and their
relation to random matrix theory.
Multiple orthogonal polynomials are a generalization of
orthogonal polynomials that have their origins in approximation
theory (Hermite-Pad\'e approximation), see e.g.~\cite{Apt2,NS}.

They enter the theory of random matrices via a generalization
of \eqref{unitarypdf} which we call a multiple orthogonal polynomial (MOP)
ensemble \cite{Kui}. We present a number of models that
give rise to a MOP ensemble, namely the model of non-intersecting
Brownian motions, the random matrix model with external source
and the two matrix model.

The MOPs are described by a Riemann-Hilbert problem that may
be used for asymptotic analysis as $n \to \infty$ by  extending
the Deift-Zhou method of steepest descent \cite{DZ}. The
extensions are non-trivial and involve either an a priori knowledge
of an underlying Riemann surface (the spectral curve) or the formulation
of a relevant equilibrium problem from logarithmic potential theory \cite{ST},
which asks for a generalization of the weighted
energy functional \eqref{weightedenergy}.

The latter approach has been succesfully applied to the
random matrix model with external source and to the two matrix model,
but only in very special cases, as will be discussed at
the end of the paper.

\section{Multiple orthogonal polynomials}

\subsection{MOP ensemble}
We will describe here multiple orthogonal polynomials
of type II, which we simply call multiple orthogonal polynomials.
There is also a dual notion of type I multiple orthogonal polynomials.

Suppose we have a finite number of weight
functions $w_1, \ldots, w_r$ on $\mathbb R$ and a multi-index
$\vec{n} = (n_1, \ldots, n_r) \in \mathbb N^r$. Associated
with these data is the monic polynomial $P_{\vec{n}}$ of degree
$|\vec{n}| = n_1 + \cdots + n_r$ so that
\begin{equation} \label{MOPintegrals}
    \int_{-\infty}^{\infty} P_{\vec{n}}(x) x^j w_k(x) \, dx = 0, \quad
    \text{ for } j=0, \ldots, n_k-1, \quad k = 1, \ldots, r.
    \end{equation}
The linear system of equations \eqref{MOPintegrals} may not be always
uniquely solvable, but in many important cases it is.
If $P_{\vec{n}}$ uniquely exists then  it is called the multiple orthogonal
polynomial (MOP) associated with the weights $w_1, \ldots, w_r$
and multi-index $\vec{n}$.

Existence and uniqueness does hold in the following situation.
Assume that
\begin{equation} \label{MOPpdf}
    \frac{1}{Z_n}
    \, \det\left[ f_j(x_k) \right]_{j,k=1, \ldots, n}
    \, \left[ \prod_{1 \leq j < k \leq n} (x_k - x_j) \right]
    \end{equation}
is a p.d.f.\ on $\mathbb R^n$, where $n = |\vec{n}|$ and
the linear span of the functions $f_1, \ldots, f_n$ is the
same as the linear span of the set of functions
\[ \{ x^{j} w_k(x) \mid j=0, \ldots, n_k-1, \, k = 1, \ldots, r \}. \]
So the assumption is that \eqref{MOPpdf} is non-negative for every choice
of $x_1, \ldots, x_n \in \mathbb R^n$, and that the normalization
constant $Z_n$ can be taken so that the integral \eqref{MOPpdf} over $\mathbb R^n$
is equal to one.
Then the MOP satisfying \eqref{MOPintegrals} exists and is given by
\[ P_{\vec{n}}(x) = \mathbb E \left[ \prod_{j=1}^n (x-x_j) \right]. \]

We call a p.d.f.\ on $\mathbb R^n$ of the form \eqref{MOPpdf}
a MOP ensemble, see \cite{Kui}.

\subsection{Correlation kernel and RH problem}
The MOP ensemble \eqref{MOPpdf} is a determinantal point process \cite{Sos}
(more precisely a biorthogonal ensemble \cite{Bor})
with a correlation kernel $K_n$ that is constructed out of multiple
orthogonal polynomials of type II and type I. It is conveniently described
in terms of the solution of a Riemann-Hilbert (RH) problem.
This RH problem for MOPs \cite{VAGK} is a generalization of the RH problem
for orthogonal polynomials due to  Fokas, Its, and Kitaev \cite{FIK}.

The RH problem asks for an $(r+1) \times (r+1)$ matrix valued function $Y$
so that
\begin{equation} \label{MOPRHproblem}
    \left\{\begin{array}{ll}
    \bullet
    {\text{ $Y : \mathbb C \setminus \mathbb R \to \mathbb C^{(r+1)\times (r+1)}$
    is analytic}}, \\
    \bullet
    {\text{ $Y$ has limiting values on $\mathbb R$, denoted by $Y_+$
    and $Y_-$, where $Y_{\pm}(x)$ is}} \\
    \quad {\text{the limit of $Y(z)$ as $z \to x \in \mathbb R$ with $\pm \Im z > 0$, satisfying }} \\
    \qquad Y_+(x) = Y_-(x) \begin{pmatrix} 1 & w_1(x) & \cdots & w_r(x) \\
        0 & 1 & \cdots & 0 \\
        \vdots & \vdots & & \vdots \\
        0 & 0 & \cdots & 1 \end{pmatrix} \text{ for } x \in \mathbb R, \\
    \bullet \ Y(z) = (I + O(1/z)) \diag \begin{pmatrix} z^n & z^{-n_1} & \cdots & z^{-n_r} \end{pmatrix}
        \text{ as } z \to \infty.
    \end{array} \right.
    \end{equation}

If the MOP $P_{\vec{n}}$ with weights $w_1, \ldots, w_r$ and multi-index
$\vec{n} = (n_1, \ldots, n_r)$ exists then the RH problem \eqref{MOPRHproblem}
has a unique solution. If the MOPs with multi-indices $\vec{n}- \vec{e}_j$
also exist, where $\vec{e}_j$ is the $j$th unit vector of length $r$,
then the first column of $Y$ consists of
\begin{equation} \label{Yfirstcolumn}
    Y_{1,1}(z) = P_{\vec{n}}(z), \qquad
    Y_{j+1,1}(z) = c_{j,\vec{n}} P_{\vec{n}- \vec{e}_j}(z), \qquad j = 1, \ldots, r
    \end{equation}
where $c_{j,\vec{n}}$ is the constant
 \[ c_{j, \vec{n}} = - 2\pi i \left[ \int_{-\infty}^{\infty}
    P_{\vec{n}-\vec{e}_j}(x) x^{n_j-1} w_j(x) dx \right]^{-1} \neq 0. \]
 The other columns of $Y$ contain Cauchy transforms
 \[ Y_{j,k+1}(z) = \frac{1}{2\pi i} \int_{-\infty}^{\infty} \frac{Y_{j,1}(x) w_k(x)}{x-z} dx,
    \qquad j=1, \ldots, r+1, \quad k = 1, \ldots, r. \]

It is a remarkable fact that the correlation kernel of the MOP ensemble \eqref{MOPpdf}
is expressed as follows in terms of the solution of the RH problem, see \cite{DaK1},
\begin{equation} \label{MOPkernel}
    K_n(x,y) = \frac{1}{2\pi i (x-y)}
    \begin{pmatrix} 0 & w_1(y) & \cdots & w_r(y) \end{pmatrix}
        Y_+^{-1}(y) Y_+(x) \begin{pmatrix} 1 \\ 0 \\ \vdots \\ 0 \end{pmatrix},
        \qquad x, y \in \mathbb R.
        \end{equation}
The inverse matrix $Y^{-1}$ contains MOPs of type I, and the
formula \eqref{MOPkernel} is essentially the Christoffel-Darboux
formula for multiple orthogonal polynomials.

Besides giving a concise formula for the correlation kernel, the
expression \eqref{MOPkernel} for the kernel gives also a possible way
to do asymptotic analysis in view of the Deift-Zhou method of steepest
descent for RH problems.

\section{Non intersecting path ensembles}

A rich source of examples of determinantal point processes is
provided by non-intersecting path ensembles. In special
cases these reduce to MOP ensembles.

\subsection{Non-intersecting Brownian motion}
Consider a one-dimensional strong Markov process with transition probability
densities $p_t(x,y)$ for $t > 0$. Suppose $n$ independent copies are given with respective
starting values $a_1 < a_2 < \cdots < a_n$ at time $t=0$
and prescribed ending values $b_1 < b_2 < \cdots < b_n$ at time
$t = T > 0$ that are conditioned not to intersect in the
full time interval $0 < t < T$.
Then by an application of a theorem of Karlin and McGregor \cite{KaMc},
the positions of the paths at an intermediate time $t \in (0,T)$ have
the joint p.d.f.
\begin{equation} \label{KarlinMcGregor}
    \frac{1}{Z_n} \, \det \left[ p_t(a_j,x_k) \right]_{1 \leq j,k\leq n} \cdot
    \det \left[p_{T-t}(x_k, b_l) \right]_{1 \leq k,l \leq n}.
    \end{equation}
In a discrete combinatorial setting the result of Karlin and McGregor is
known as the Lindstrom-Gessel-Viennot theorem.

\begin{figure}[t]
\centering
\begin{overpic}[width=8cm]{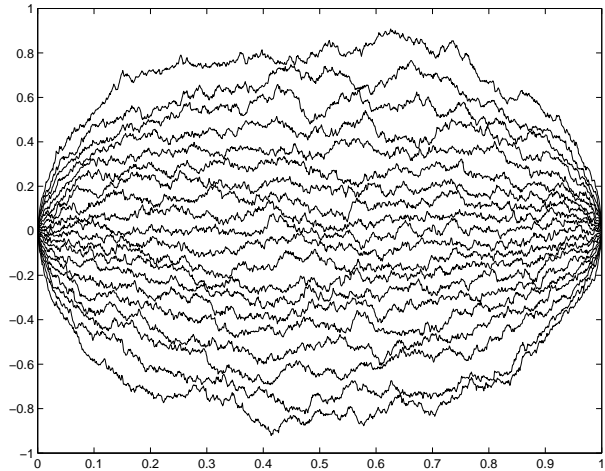}
\end{overpic}
\caption{Non-intersecting Brownian bridges starting and ending at $0$.
At any intermediate time $t \in (0,1)$ the positions of the paths have the same
distribution as the (appropriately rescaled) eigenvalues of an $n \times n$ GUE matrix. \label{Figure1}}
\end{figure}

The density function \eqref{KarlinMcGregor} is a biorthogonal ensemble,
which in very special cases reduces to the form \eqref{MOPpdf} of
a MOP ensemble.

An example is the case of Brownian motion (actually Brownian bridges)
with the transition probability density
\[ p_t(x,y) = \frac{1}{\sqrt{2\pi t} \sigma} e^{-\frac{(x-y)^2}{2t}},
        \qquad t > 0. \]
In the confluent limit where all $a_j \to 0$ and all $b_l \to 0$
the p.d.f.\  \eqref{KarlinMcGregor}  turns into
\[ \frac{1}{Z_n} \prod_{1 \leq j < k \leq n} (x_k-x_j)^2 \prod_{j=1}^n e^{- \frac{T}{2t(T-t)} x_j^2} \]
with a different constant $Z_n$. This is up to trivial scaling
the same as the p.d.f.\ \eqref{GUEdensity} for the eigenvalues of an $n \times n$ GUE matrix.

If however, we let all $b_l \to 0$ and choose
only $r$ different starting values, denoted
by $a_1, \ldots, a_r$, and $n_j$ paths start at $a_j$,
then \eqref{KarlinMcGregor} turns into a MOP ensemble with weights
\begin{equation} \label{multipleHermite}
    w_j(x) = e^{- \frac{T}{2t(T-t)} x^2 + \frac{a_j}{t} x}, \qquad j=1,\ldots, r,
    \end{equation}
and multi-index $(n_1, \ldots, n_r)$.
This is a multiple Hermite ensemble, since the associated MOPs
are multiple Hermite polynomials \cite{ABV}

\begin{figure}[t]
\centering
\begin{overpic}[width=8cm]{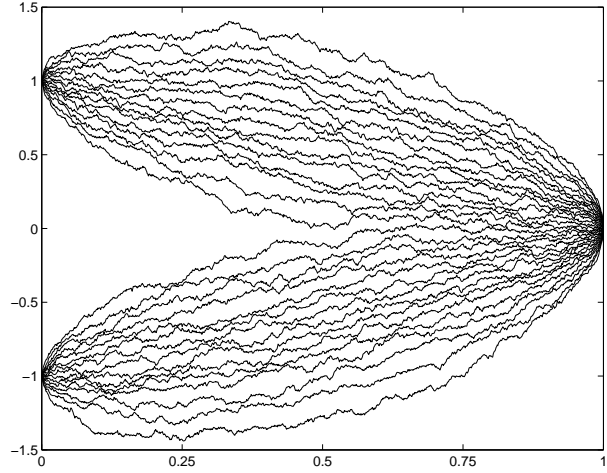}
\end{overpic}
\caption{Non-intersecting Brownian bridges starting at two different values
and ending at $0$.
At any time $t \in (0,1)$ the positions of the paths have the same
distribution as the eigenvalues of an $n \times n$ GUE matrix with external source.
The distribution is a multiple Hermite ensemble with two Gaussian weights
\eqref{multipleHermite}. \label{Figure2}}
\end{figure}

\subsection{Non-intersecting squared Bessel paths}
The squared Bessel process is another one-dimensional Markov process
which gives rise to a MOP ensemble. The squared Bessel process is
a Markov process on $[0,\infty)$, depending on a parameter $\alpha > -1$,
with transition probability density
\[ p_t(x,y) = \frac{1}{2t} \left( \frac{y}{x} \right)^{\alpha/2}
    e^{- \frac{1}{2t} (x+y)} I_{\alpha} \left( \frac{\sqrt{xy}}{t} \right),
        \qquad x, y > 0,   \]
where $I_{\alpha}$ is the modified Bessel function of first kind of order $\alpha$.
In the limit where all $a_j \to a > 0$ and $b_j \to 0$ the
p.d.f.\ \eqref{KarlinMcGregor} for the positions of the paths
at time $t \in (0,T)$ is a MOP ensemble with two weights
\begin{align*}
    w_1(x) & = x^{\alpha/2} e^{- \frac{T}{2t(T-t)} x} I_{\alpha} \left( \frac{\sqrt{ax}}{t} \right)  \\
    w_2(x) & = x^{(\alpha+1)/2} e^{- \frac{T}{2t(T-t)} x} I_{\alpha+1} \left( \frac{\sqrt{ax}}{t} \right)
    \end{align*}
and multi-index $(n_1, n_2)$ where $n_1 = \lceil n/2  \rceil$
and $n_2 = \lfloor n/2 \rfloor$, see \cite{KMW}. In the limit $a \to 0$ this
further reduces to an orthogonal polynomial ensemble for a Laguerre weight.

\section{Random matrix models}

The random matrix model with external source, and the
two matrix model also give rise to MOP ensembles.

\subsection{Random matrices with external source}

The Hermitian matrix model with external source is
the probability measure
\begin{equation} \label{sourcemodel}
    \frac{1}{Z_n}  e^{-n \Tr (V(M) - AM)} dM
    \end{equation}
on $n \times n$ Hermitian matrices, where the external source
$A$ is a given Hermitian $n \times n$ matrix.
This is a modification of the usual Hermitian matrix
model, in which the unitary invariance is broken \cite{BH}, \cite{ZJ}.

Due to the Harish-Chandra/Itzykson-Zuber integral formula \cite{IZ},
it is possible to integrate out the eigenvectors explicitly.
In case the eigenvalues $a_1, \ldots, a_n$ of $A$ are all
distinct, we obtain the explicit p.d.f.\
\[ \frac{1}{Z_n} \det \left[ e^{n a_i x_j} \right]_{1 \leq i, j \leq n} \cdot
    \prod_{1 \leq j< k \leq n} (x_k-x_j) \cdot \prod_{j=1}^n e^{-n V(x_j)} \]
for the eigenvalues of $M$. In case that $a_1, \ldots, a_r$ are the distinct eigenvalues
of  $A$, with respective multiplicities
$n_1, \ldots, n_r$, then the eigenvalues of $M$ are distributed as
a MOP ensemble \eqref{MOPpdf} with weights
\begin{equation} \label{sourceweights}
    w_j(x) =  e^{-n (V(x) - a_j x)}, \qquad j=1, \ldots, r
    \end{equation}
and multi-index $(n_1, \ldots, n_r)$, see \cite{BK0}

For the case $V(x) = \frac{1}{2} x^2$ the external source model \eqref{sourcemodel}
is equivalent to the model of non-intersecting Brownian motions with
several starting points and one ending point, cf.\ \eqref{multipleHermite}.

\subsection{Two matrix model}
The Hermitian two matrix model
\begin{equation} \label{twomatrixmodel}
    \frac{1}{Z_n} e^{-n \Tr(V(M_1) + W(M_2) - \tau M_1 M_2)} \, dM_1 dM_2
    \end{equation}
is a probability measure defined on couples $(M_1, M_2)$ of $n \times n$
Hermitian matrices. Here $V$ and $W$ are two potentials (typically polynomials)
and $\tau \neq 0$ is a coupling constant. The model is of great interest
in $2$d quantum gravity \cite{DKK}, \cite{IZ}, \cite{Kaz},
as it allows for a large class of critical phenomena.

The eigenvalues of the matrices $M_1$ and $M_2$ are fully described
by biorthogonal polynomials.
These are two sequences $(P_{k,n})_k$ and $(Q_{j,n})_j$
of monic polynomials, $\deg P_{k,n} = k$,
$\deg Q_{j,n} = j$, such that
\begin{equation} \label{biorthogonal}
    \int_{-\infty}^{\infty} \int_{-\infty}^{\infty} P_{k,n}(x) Q_{j,n}(y)
    e^{-n(V(x) + W(y) - \tau xy)} \, dx \, dy = 0,
    \qquad \text{if } j \neq k,
    \end{equation}
see e.g.\ \cite{BEH}, \cite{EMc}, \cite{Meh}, \cite{MS}.

If $W$ is a polynomial then the biorthogonality conditions \eqref{biorthogonal}
can be seen as multiple orthogonal polynomial conditions with
respect to $r = \deg W -1$ weights
\begin{equation} \label{twomatrixweights}
    w_{j,n}(x) = e^{-n V(x)} \int_{-\infty}^{\infty} y^j e^{-n (W(y) - \tau xy)} dy,
    \qquad j = 0, \ldots, r-1,
    \end{equation}
see \cite{KMc}. Furthermore, the eigenvalues of $M_1$ are a MOP ensemble \eqref{MOPpdf}
with the weights \eqref{twomatrixweights} and multi-index $\vec{n} = (n_0, \ldots, n_{r-1})$
with $n_j = \lceil n/r \rceil$ for $j=0, \ldots, q-1$
and $n_j = \lfloor n/r \rfloor$ for $j= q, \ldots, r-1$
if $n = pr + q$ with $p$ and $0 \leq  q < r$ non-negative integers, see \cite{DuK2}
for the case where $W(y) = \frac{y^4}{4}$.

\section{Large $n$ behavior and critical phenomena}

We discuss the large $n$ behavior in the above described models.

\subsection{Non-intersecting Brownian motion}

In order to have interesting limit behavior as $n \to \infty$ in
the non-intersecting Brownian motion model we scale
the time variables $T \mapsto 1/n$, $t \mapsto t/n$, so that $0 < t < 1$.
In the case of one starting value and one ending value, see
Figure \ref{Figure1}, the paths will fill out an ellipse as $n \to \infty$.

\begin{figure}[t]
\centering
\begin{overpic}[width=8cm]{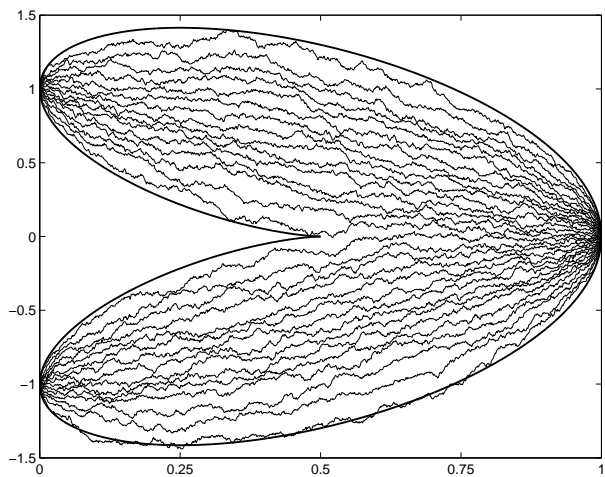}
\end{overpic}
\caption{Non-intersecting Brownian bridges starting at two different values
and ending at $0$. As $n \to \infty$, the paths fill out a heart-shaped
domain. Critical behavior at the cusp point is desribed by the Pearcey
kernel \eqref{Pearceykernel}.  \label{Figure3}}
\end{figure}

In the situation of Figure \ref{Figure2} the paths fill out
a heart-shaped region as $n \to \infty$, as shown in Figure \ref{Figure3}. New critical behavior
appears at the cusp point where the two groups of paths
come together and merge into one.

Around the critical time the correlation kernels have
a double scaling limit, which is given by the one-parameter family of Pearcey kernels
\begin{equation} \label{Pearceykernel}
    \mathcal P(x,y; b) = \frac{p(x) q''(y) - p'(x) q'(y) + p''(x) q(y) - b p(x) q(y)}{x-y},
        \qquad b \in \mathbb R,
    \end{equation}
where $p$ and $q$ are solutions of the Pearcey equations
$p'''(x) = x p(x) - b p'(x)$ and $q'''(y) = yq(y) + b q'(y)$.
This kernel was first identified by Br\'ezin and Hikami \cite{BH}
who also gave the double integral representation
\begin{equation} \label{Pearceydoubleintegral}
    \mathcal P(x,y) =
    \frac{1}{(2\pi i)^2} \int_C \int_{-i\infty}^{i\infty}
        e^{- \frac{1}{4} s^4 + \frac{b}{2} s^2 - ys + \frac{1}{4} t^4 - \frac{b}{2} t^2 +x t}
            \frac{ds \, dt}{s-t}
            \end{equation}
where the contour $C$ consists of the rays from $\pm \infty e^{i\pi/4}$ to $0$ and
the rays from $0$ to $\pm \infty e^{-i \pi/4}$.

Consideration of multiple times near the critical time leads to an extended Pearcey kernel
and the Pearcey process given by Tracy and Widom \cite{TW2}.

As already noted above, the model of non-intersecting Brownian motion with
two starting points and one ending point is related to the
Gaussian random matrix model with external source
\begin{equation} \label{GUEsourcemodel}
    \frac{1}{Z_n} e^{-n \Tr (\frac{1}{2} M^2 - AM)} \, dM,
\end{equation}
with external source
\begin{equation} \label{externalsourceA}
     A = \diag(\underbrace{a, \ldots, a}_{n/2 \text{ times}},
        \underbrace{-a, \ldots, -a}_{n/2 \text{times}}).
        \end{equation}
In this setting the critical $a$-value is $a_{crit}=1$
and the Pearcey kernel \eqref{Pearceykernel} arises as $n \to \infty$ with
$a = 1 + \frac{b}{2\sqrt{n}}$. In \cite{BK3} this was studied with the
use of the Riemann-Hilbert problem \eqref{MOPRHproblem} for multiple
Hermite polynomials with two weights $e^{- n( \frac{1}{2} x^2 \pm ax)}$.
The asymptotic analysis as $n \to \infty$ was done with an extension of
the Deift-Zhou method of steepest descent \cite{DZ}
to the case of a $3 \times 3$ matrix valued RH problem.
See also \cite{BK1} and \cite{ABK} for a steepest descent analysis of
the RH problem in the non-critical regimes $a > 1$ and $0 < a < 1$, respectively.

Another interesting asymptotic regime is the model of non-intersecting
Brownian motion with outliers. In this model a rational  modification of the Airy kernel
appears that was first described in \cite{BBAP} in the context
of complex sample covariance matrices, see also \cite{ADvM}.

\subsection{Random matrices with external source}

If $V$ is quadratic in the random matrix model with external source
\eqref{sourcemodel} then this model can be mapped to the model
of non-intersecting Brownian motions. Progress on this model
beyond the quadratic case is due to McLaughlin \cite{McL}
who found the spectral curve for the quartic potential $V(x) = \frac{1}{4}x^4$
and for $a$ sufficiently large (again $A$ is as in \eqref{externalsourceA}).

A method based on a vector equilibrium problem was introduced
recently by Bleher, Delvaux and Kuijlaars \cite{BDK}.
The vector equilibrium problem extends the equilibrium problem
for the weighted energy \eqref{weightedenergy}
that is important for the unitary ensembles and which is crucial
in the steepest descent analysis of the RH problem for orthogonal
polynomials \cite{DKMVZ}.

In \cite{BDK} it is assumed that $V$ is an even polynomial,
and that $A$ is again given as in \eqref{externalsourceA}.
The vector equilibrium problem involves two measures
$\mu_1$ and $\mu_2$, and it asks to minimize
the energy functional
\begin{multline} \label{externalsourcevectorEP}
    \iint \log \frac{1}{|x-y|} d\mu_1(x) d\mu_1(y) +
    \iint \log \frac{1}{|x-y|} d\mu_2(x) d\mu_2(y) \\
    -
    \iint \log \frac{1}{|x-y|} d\mu_1(x) d\mu_2(y)
    + \int (V(x) - a |x|) \, d\mu_1(x)
    \end{multline}
where $\mu_1$ is on $\mathbb R$ with $\int d\mu_1 = 1$,
$\mu_2$ is on $i \mathbb R$ (the imaginary axis) with $\int d\mu_2 = 1/2$,
and in addition $\mu_2 \leq \sigma$, where $\sigma$ is the
measure on $i \mathbb R$  with constant density
\begin{equation} \label{externalsourceconstraint}
    \frac{d\sigma}{|dz|} = \frac{a}{\pi}.
    \end{equation}

There is a unique minimizer, and the density $\rho_1$ of
the measure $\mu_1$ is the limiting mean eigenvalue density
\[ \rho_1(x) = \lim_{n \to \infty} \frac{1}{n} K_n(x,x) \]
where $K_n$ is the correlation kernel of the MOP ensemble with
weights $e^{-n (V(x) \pm ax)}$. The RH problem \eqref{MOPRHproblem}
is analyzed in the large $n$ limit with the Deift/Zhou steepest descent
method in which the minimizers from the vector equilibrium
problem play a crucial role.

The upper constraint $\mu_2 \leq \sigma$ is not active for large enough $a$
and in that case the support of $\mu_1$ has a gap around $0$.
For smaller values of $a$ the constraint $\sigma$ is active along
an interval $[-ic, ic]$, $c > 0$, on the imaginary axis.
Critical phenomena take place when either the constraint becomes
active, or the gap around $0$ closes, or both. If one
of these two phenomena happens, then this generically will be
a phase transition of the Painlev\'e II type that was
described in the unitary matrix model in \cite{BI2}
and \cite{CK}. If the two phenomena happen simultaneously
then this is expected to be phase transition of the Pearcey type
which, if true, would be a confirmation of the universality of
the Pearcey kernels \eqref{Pearceykernel} at the closing of a gap
\cite{BH}.

Both kinds of transitions are valid in the external source model with even
quartic potential $V(x) = \frac{1}{4} x^4 - \frac{t}{2} x^2$, see \cite{BDK}.
For the particular value $t = \sqrt{3}$, there is a passage from the Painlev\'e II
transition (for $t > \sqrt{3}$) to the Pearcey transition (for $t < \sqrt{3}$).
The description of the phase transition for $t = \sqrt{3}$ remains open.

\subsection{Two matrix model}
In \cite{DuK2} Duits and Kuijlaars applied the steepest descent analysis
to the RH problem for the two matrix model \eqref{twomatrixmodel} with quartic potential
\begin{equation} \label{quarticW}
    W(y) = \frac{1}{4} y^4
    \end{equation}
and for $V$ an even polynomial. The corresponding MOP ensemble has
three weights of the form \eqref{twomatrixweights} and the RH problem
\eqref{MOPRHproblem} is of size $4 \times 4$.
Again a vector equilibrium problem plays a crucial  role.

The vector equilibrium problem in \cite{DuK2} involves three
measures $\mu_1$, $\mu_2$ and $\mu_3$. It asks to minimize
the energy functional
\begin{multline} \label{twomatrixmodelvectorEP}
   \sum_{j=1}^3 \iint \log \frac{1}{|x-y|} d\mu_j(x) d\mu_j(y) \\
   - \sum_{j=1}^2 \iint \log \frac{1}{|x-y|} d\mu_j(x) d\mu_{j+1}(y)
     + \int (V(x) -  \frac{3}{4} |\tau x|^{4/3}) \, d\mu_1(x)
    \end{multline}
among measures $\mu_1$ on $\mathbb R$ with $\int d\mu_1 = 1$,
$\mu_2$ on $i \mathbb R$ with $\int d\mu_2 = 2/3$
and $\mu_3$ on $\mathbb R$ with $\int d\mu_3 = 1/3$.
In addition $\mu_2 \leq \sigma$
where $\sigma$ is a given measure on $i\mathbb R$ with density
\begin{equation} \label{twomatrixmodelconstraint}
    \frac{d\sigma}{|dz|} = \frac{\sqrt{3}}{2\pi} |\tau|^{4/3} |z|^{1/3},
    \qquad z \in i \mathbb R.
    \end{equation}

There is a unique minimizer and the density $\rho_1$ of the first
measure $\mu_1$ is equal to the limiting mean eigenvalue density
of the matrix $M_1$ in the two matrix model. In addition, the usual
scaling limits (sine kernel and Airy kernel) are valid in
the local eigenvalue regime, see \cite{DuK2}. However there
is no new critical behavior in the two matrix model
with $W$ is given by \eqref{quarticW}.

New multicritical behavior is predicted in \cite{DKK}
for more general potentials.
For the more general quartic potential $W(y) = \frac{1}{4} y^4 - \frac{t}{2} y^2$
an approach based on a modification of
the vector equilibrium problem \eqref{twomatrixmodelvectorEP}
is under current investigation.

\bibliographystyle{amsalpha}

\end{document}